\documentclass[11pt]{article}
\usepackage[margin=1in]{geometry}
\usepackage[T1]{fontenc}
\usepackage[utf8]{inputenc}
\usepackage{amsmath,amssymb,amsthm,array,booktabs,microtype,enumitem}
\usepackage{cite}
\usepackage{tikz}
\usetikzlibrary{arrows.meta,decorations.markings}
\usepackage[colorlinks=true,urlcolor=blue,linkcolor=black,citecolor=black]{hyperref}
\hypersetup{pdftitle={The Cranston\textendash Li\textendash Su\textendash Wang\textendash Xu (2026) conjecture is true},pdfauthor={Dave Goel},
  pdfkeywords={nowhere-zero flow, group flow, reconfiguration, flow-connected graph, Eulerian graph}}

\newtheorem{theorem}{Theorem}[section]
\newtheorem{lemma}[theorem]{Lemma}
\newtheorem{corollary}[theorem]{Corollary}
\newtheorem{proposition}[theorem]{Proposition}
\theoremstyle{definition}
\newtheorem{definition}[theorem]{Definition}
\newtheorem{remark}[theorem]{Remark}
\newcommand{\Z}{\mathbb{Z}}
\newcommand{\cF}{\mathcal{F}}
\newcommand{\cC}{\mathcal{C}}
\newcommand{\one}{\mathbf{1}}
\tikzset{
  vtx/.style={circle,fill=black,inner sep=1.4pt},
  arr/.style={-{Stealth[length=2.2mm]},thick},
  stub/.style={thick},
  cell/.style={draw,minimum width=8mm,minimum height=6.4mm,text height=1.5ex,text depth=0.35ex,inner sep=0pt,font=\scriptsize},
  midarr/.style={thick,postaction={decorate},decoration={markings,mark=at position 0.62 with {\arrow{Stealth[length=2.2mm]}}}},
  yes/.style={cell,fill=green!35},
  open/.style={cell,fill=orange!25},
  no/.style={cell,fill=black!8,very thick},
}
\makeatletter\@secpenalty=0\makeatother
\title{The Cranston--Li--Su--Wang--Xu (2026) conjecture is true}
\author{Dave Goel\thanks{\raggedright All theorems of this paper are proved in Lean~4 from Lean's three standard
axioms alone, with no other result encoded as an axiom, and checked with the Lean comparator:
\url{https://github.com/david44hello/papers_supporting_materials/tree/cranston-li-su-wang-xu-2026-conjecture-v2.00/cranston-li-su-wang-xu-2026-conjecture}}}
\date{October 2026\\[0.4ex] {\normalsize Revised: November 10, 2026}}

\begin{document}
\maketitle\vspace*{-1em}{\footnotesize\centering\itshape This paper has an expanded version, with further results~\cite{G}: \textup{\url{https://doi.org/10.5281/zenodo.23605838}}\par}

\begin{abstract}
Let $A$ be a finite abelian group. A nowhere-zero $A$-flow of a graph can be changed along a cycle: going around the
cycle, add a fixed non-zero element of $A$ to the value of
every edge that points forward and subtract it from the value of every edge that points backward. The Cranston--Li--Su--Wang--Xu (2026) conjecture on the reconfiguration of nowhere-zero flows states that if $A$ has at least four elements, then in every
graph in which all degrees are even, any two nowhere-zero
$A$-flows are joined by a chain of such changes through nowhere-zero $A$-flows. It also states that, for every integer $k\ge 4$, the same holds for integer flows whose
values are between $1$ and $k-1$ in absolute value, where each change adds a non-zero integer. We prove this conjecture. The proof is by induction on the number of edges. If a flow
enters a vertex along one edge and leaves it along another edge with the same value, the two edges can be
replaced by a single edge; this gives a smaller graph. The main step is to change every flow into a flow with such a pair of edges.
For $A$-flows we do this when all degrees are at most $4$; repeatedly replacing a vertex of degree at least $6$ by a
triangle, following Cranston et al., reduces the general case to this one.
For integer flows, we show that, in a graph in which all degrees are even, any two nowhere-zero flows with values in the cyclic group $\Z_k$, $k\ge 4$, are joined by
changes that add only $1$, $-1$, $2$ or $-2$, and that such changes lift to integer flows. All theorems of this paper
are also proved in the Lean~4 proof assistant.
\end{abstract}


\section{Introduction}\label{sec:intro}

This paper proves the Cranston--Li--Su--Wang--Xu (2026) conjecture \cite[Conjecture~4.1]{CLSWX} on the
reconfiguration of nowhere-zero flows in graphs in which every vertex has even degree, that is, on whether every
nowhere-zero flow can be turned into any other one by changes along cycles that keep the flow nowhere-zero.

Let $G$ be a finite graph; loops and parallel edges are allowed. For a finite abelian group $A$ (for example the cyclic group $\Z_m$ of the integers modulo $m$), a nowhere-zero
$A$-flow of $G$ assigns to every edge, in a fixed reference direction, a non-zero element of $A$, so that at every
vertex the sum of the values on the entering edges equals the sum on the leaving edges. Two such flows differ by a
\emph{step} with \emph{value} $t\ne 0$ if there is a cycle, traversed in one of its two directions, such that their
difference is $t$ on every edge of the cycle whose reference direction agrees with the traversal, $-t$ on every other
edge of the cycle, and $0$ on every edge not on the cycle. The graph $\cF(G,A)$ has the nowhere-zero $A$-flows of $G$ as its vertices, two of them being
adjacent if they differ by a step, and $G$ is \emph{$A$-flow-connected} if $\cF(G,A)$ is connected. For an integer
$k\ge 2$, the graph $\cF(G,k)$ is defined in the same way with integers. Its vertices are the \emph{$k$-flows} of $G$, that is, the integer flows $F$ of $G$
with $1\le|F(e)|\le k-1$ on every edge $e$ (they are usually called nowhere-zero $k$-flows), and two of them are adjacent if they differ by a
step whose value is a non-zero integer. The graph $G$ is \emph{$k$-flow-connected} if $\cF(G,k)$ is connected.
Figure~\ref{fig:step} shows a step.

\begin{figure}[ht]
\centering
\begin{tikzpicture}[scale=1.15]
  \foreach \side/\dx in {L/0,R/5.6} {
    \node[vtx] (a\side) at (\dx+0,0) {};
    \node[vtx] (b\side) at (\dx+2,0) {};
    \node[vtx] (c\side) at (\dx+2,2) {};
    \node[vtx] (d\side) at (\dx+0,2) {};
    \draw[->,thin] (\dx+1,1) ++(-150:0.45) arc (-150:150:0.45);
    \node[font=\scriptsize] at (\dx+1,1) {$Q$};
  }
  \draw[arr] (aL) -- node[below,font=\small] {$f(e_1)$} (bL);
  \draw[arr] (bL) -- node[right,font=\small] {$f(e_2)$} (cL);
  \draw[arr] (dL) -- node[above,font=\small] {$f(e_3)$} (cL);
  \draw[arr] (aL) -- node[left,font=\small] {$f(e_4)$} (dL);
  \draw[arr] (aR) -- node[below,font=\small] {$f(e_1)+t$} (bR);
  \draw[arr] (bR) -- node[right,font=\small] {$f(e_2)+t$} (cR);
  \draw[arr] (dR) -- node[above,font=\small] {$f(e_3)-t$} (cR);
  \draw[arr] (aR) -- node[left,font=\small] {$f(e_4)-t$} (dR);
  \node[font=\small] at (1,-0.9) {the flow $f$};
  \node[font=\small] at (6.6,-0.9) {the flow after the step};
  \draw[-{Stealth[length=3mm]},very thick] (3.15,1) -- node[above,font=\small]{step} (4.15,1);
\end{tikzpicture}
\caption{A step with value $t$ along a cycle $Q$ with four edges, traversed counter-clockwise. The arrow on each edge
is its reference direction. On the two edges whose reference direction agrees with the traversal, $t$ is added to the
value; on the other two, $t$ is subtracted; the edges not on $Q$ (not drawn) keep their values. The flow condition
still holds at every vertex.}
\label{fig:step}
\end{figure}

Esperet, Hendrey, Lagoutte, Marseloo, Norin and Steiner \cite{EHLMNS} introduced the graphs $\cF(G,A)$ and $\cF(G,k)$
and asked for which graphs $G$, groups $A$ and integers $k$ they are connected. Cranston, Li, Su, Wang and Xu
\cite{CLSWX} call a graph in which all degrees are even Eulerian. Since ``Eulerian'' often also requires the graph to
be connected, we call a graph \emph{even} if every vertex has even degree; an even graph need not be connected.

\begin{theorem}\label{thm:main}
Let $G$ be a graph in which every vertex has even degree.
\begin{itemize}
\item[(a)] $G$ is $A$-flow-connected for every finite abelian group $A$ with $|A|\ge 4$.
\item[(b)] $G$ is $k$-flow-connected for every integer $k\ge 4$.
\end{itemize}
\end{theorem}

Theorem~\ref{thm:main} is \cite[Conjecture~4.1]{CLSWX}. Part~(a) is Theorem~\ref{thm:all}
(Section~\ref{sec:proofs}), and part~(b) is Theorem~\ref{thm:int} (Section~\ref{sec:int}). The table
and Figure~\ref{fig:result} compare Theorem~\ref{thm:main} with earlier results.

\begin{center}
\small
\begin{tabular}{@{}>{\raggedright\arraybackslash}p{0.68\textwidth}@{\hspace{0.03\textwidth}}>{\raggedright\arraybackslash}p{0.29\textwidth}@{}}
\toprule
Best result before this paper & New result\\
\midrule
\multicolumn{2}{@{}l}{\emph{Which graphs are $A$-flow-connected, and for which finite abelian groups $A$?}}\\
\addlinespace[0.3ex]
All even graphs: even $|A|$ \cite[Theorem~4.3]{CLSWX}; $|A|\ge 5{,}251{,}339$ \cite[Corollary~4.5]{CLSWX}; for
$|A|\ge 5$, reduced to $4$-regular graphs \cite[Observation~4.6]{CLSWX}.\newline
All graphs without a bridge: $|A|\ge 3024$ \cite[Theorem~1.1]{GoelGroups}; some groups of smaller order, including
$\Z_{55}$ \cite[Theorem~4.5]{GoelEarlier}.\newline
Even loopless planar graphs: all $|A|\ge 4$, by \cite{CLSWX} and \cite{EHLMNS} (see below).
& All even graphs: all $|A|\ge 4$ (Theorem~\ref{thm:main}~(a)).\\
\addlinespace
\multicolumn{2}{@{}l}{\emph{Which graphs are $k$-flow-connected, and for which integers $k$?}}\\
\addlinespace[0.3ex]
All even graphs: even $k$ \cite[Theorem~4.2]{CLSWX}.\newline
All graphs without a bridge and without loops: some odd $k$, for example $k=729$
\cite[Theorem~7.4]{GoelEarlier}.\newline
Even loopless planar graphs: all $k\ge 4$, by \cite{CLSWX} and \cite{EHLMNS} (see below).
& All even graphs: all $k\ge 4$ (Theorem~\ref{thm:main}~(b)).\\
\bottomrule
\end{tabular}
\end{center}

An even graph has no bridge, so the results in the table for graphs without a bridge, from \cite{GoelEarlier} and
\cite{GoelGroups}, apply to even graphs. Such results
go back to \cite[Theorems~6.15 and~6.23]{EHLMNS}, for the group $\Z_2^8$ and for all groups of sufficiently large order.
Neither \cite{GoelEarlier} nor \cite{GoelGroups} is used here.

For even planar graphs without loops (the graphs of \cite{EHLMNS} have no loops), the conjecture already followed from
earlier results, applied to each connected component. A
connected even plane graph $G$ without loops has no bridge, and its dual is bipartite; so the simple graph underlying
the dual is bipartite and planar, hence $3$-degenerate by Euler's formula (as noted in \cite[Section~5]{EHLMNS}). By
\cite[Corollary~5.8]{EHLMNS}, if the simple graph underlying the dual is $d$-degenerate, then $G$ is
$A$-flow-connected for $|A|\ge d+2$ and $k$-flow-connected for $k\ge d+2$; here $d=3$. The cases $|A|=4$ and $k=4$ are
\cite[Theorems~4.3 and~4.2]{CLSWX}.

\begin{figure}[ht]
\centering
\begin{tikzpicture}[x=8mm,y=8mm]
  \node[font=\small\bfseries] at (3.5,2.3) {Known before this paper};
  \node[font=\small\bfseries] at (12.7,2.3) {Known after this paper};
  \foreach \lab/\x in {3/0,4/1,5/2,6/3,7/4,8/5,9/6} {
    \node[font=\scriptsize] at (\x,1.5) {$\lab$};
    \node[font=\scriptsize] at (\x+9.2,1.5) {$\lab$};
  }
  \node[font=\scriptsize] at (7,1.5) {$\cdots$};
  \node[font=\scriptsize] at (16.2,1.5) {$\cdots$};
  \node[font=\footnotesize,anchor=east] at (-0.65,0.8) {groups, $|A|=$};
  \node[font=\footnotesize,anchor=east] at (-0.65,0) {integers, $k=$};
  \foreach \y in {0.8,0} {
    \node[no] at (0,\y) {no};
    \foreach \x in {1,3,5} \node[yes] at (\x,\y) {yes};
    \foreach \x in {2,4,6} \node[open] at (\x,\y) {open};
    \node[cell] at (7,\y) {$\cdots$};
    \node[no] at (9.2,\y) {no};
    \foreach \x in {1,...,7} \node[yes] at (\x+9.2,\y) {yes};
  }
\end{tikzpicture}
\caption{Is every even graph $A$-flow-connected (rows ``groups'', by the order $|A|$) and $k$-flow-connected (rows
``integers'', by $k$)? For $3$ the answer is no \cite[Theorem~2.1]{CLSWX}; see Figure~\ref{fig:dipole}. Left: what
was known before this paper for the values drawn \cite[Theorems~4.2 and~4.3]{CLSWX}. Right: what is known after this paper; every ``yes'' there follows from
Theorem~\ref{thm:main}.}
\label{fig:result}
\end{figure}

The number $4$ in Theorem~\ref{thm:main} cannot be replaced by $3$. By \cite[Theorem~2.1]{CLSWX}, for a graph with a
$3$-flow, the following are equivalent: it is
$\Z_3$-flow-connected; it is $3$-flow-connected; all its degrees are even and no two of its vertices are joined by more
than four edge-disjoint paths. The even graph with two vertices and six parallel edges has a
$3$-flow (Figure~\ref{fig:dipole}), and its
two vertices are joined by six edge-disjoint paths. So it is neither $\Z_3$-flow-connected nor $3$-flow-connected
(for $\Z_3$, this is also noted in \cite[Section~5]{EHLMNS}).

\begin{figure}[ht]
\centering
\begin{tikzpicture}[scale=1.0]
  \node[vtx,label=left:{$x$}] (x) at (0,0) {};
  \node[vtx,label=right:{$y$}] (y) at (5,0) {};
  \foreach \bend/\val in {66/1,40/1,14/1,-14/-1,-40/-1,-66/-1} {
    \draw[midarr] (x) to[bend left=\bend] node[pos=0.36,fill=white,inner sep=1pt,font=\scriptsize] {$\val$} (y);
  }
\end{tikzpicture}
\caption{The even graph with two vertices and six parallel edges, with a $3$-flow: all edges are directed from $x$ to
$y$; three have the value $1$ and three have the value $-1$.}
\label{fig:dipole}
\end{figure}

\paragraph{Two further theorems.} For the cyclic groups $\Z_m$ with $m\ge 4$, any two nowhere-zero $\Z_m$-flows of an even graph are joined by a
chain of steps whose values are in $\{1,-1,2,-2\}$ (Theorem~\ref{thm:cyclic}); for even $m$, $\Z_m$-flow-connectedness, with steps of any value,
was known \cite[Theorem~4.3]{CLSWX}. Part~(a) rests on Theorem~\ref{thm:deg4}: every even graph in which all degrees are at most
$4$ is $A$-flow-connected for every finite abelian group $A$. For $|A|\le 3$ this is not new: it is trivial for
$|A|\le 2$, as a graph then has at most one nowhere-zero $A$-flow, and for $A=\Z_3$ it follows from
\cite[Theorem~2.1]{CLSWX}, because an even graph has a $3$-flow (the value $1$ along an Eulerian orientation) and, when all
its degrees are at most $4$, no two of its vertices are joined by more than four edge-disjoint paths.

\paragraph{Outline of the proof.}
The proof is by induction on the number of edges, and its main tool is the balanced pair. A flow has a balanced pair
at a vertex $v$ if it enters $v$ along an edge $e$ from a vertex $u$ and leaves $v$ along another edge $e'$ to a vertex $w$ with the same value
(Definition~\ref{def:balanced}). The two edges can then be replaced by a single edge $uw$
(Figure~\ref{fig:split}), and the flows with this balanced pair correspond to the flows of the smaller
graph (Lemma~\ref{lem:split}). By the induction hypothesis the smaller graph is flow-connected, so any two flows with
a balanced pair on the same two edges are joined by a chain of steps.

Flows whose balanced pairs are on different edges
are joined through constant flows: a constant flow has the same value $a\ne 0$ on every edge, in the direction of an
Eulerian orientation (an orientation in which, at every vertex, as many edges enter as leave; one exists because every
degree is even). Any two edges at a vertex of a loopless even graph form a balanced pair in some constant flow with
value $a$, and all constant flows with value $a$ are joined to each other, for every $a\ne 0$ when steps of all values
are allowed and for $a=1$ when only steps with values $\pm1$ and $\pm2$ are allowed (Lemma~\ref{lem:const}).
So all flows with a balanced pair are joined to each other, and it remains to join every flow to one of them. This
last statement is called \textup{(LOCAL)} in Proposition~\ref{prop:red}, which states the reduction precisely.

Constant flows along Eulerian orientations and splitting off two edges at a vertex both come from \cite{CLSWX}. In the proofs of \cite[Theorems~4.2 and~4.3]{CLSWX} every flow is brought
to a constant flow along an Eulerian orientation, and any two Eulerian orientations are joined by reversals of directed
cycles \cite[Lemma~2.3]{CLSWX} (see Lemma~\ref{lem:const}~(a)). The proof of \cite[Lemma~2.3]{CLSWX}
also splits off two edges at a vertex; Lemma~\ref{lem:split} extends this step from Eulerian orientations to all
flows with a balanced pair.

\begin{figure}[ht]
\centering
\begin{tikzpicture}[scale=1.1]
  \node[vtx,label=below:{$v$}] (v) at (0,0) {};
  \node[vtx,label=left:{$u$}] (u) at (-1.8,0.9) {};
  \node[vtx,label=right:{$w$}] (w) at (1.8,0.9) {};
  \draw[arr] (u) -- node[above=2pt,font=\small] {$x$} node[below=2pt,font=\small] {$e$} (v);
  \draw[arr] (v) -- node[above=2pt,font=\small] {$x$} node[below=2pt,font=\small] {$e'$} (w);
  \draw[stub,gray] (v) -- (-1.0,-1.1); \draw[stub,gray] (v) -- (1.0,-1.1);
  \draw[-{Stealth[length=3mm]},very thick] (2.6,0.2) -- (3.6,0.2);
  \begin{scope}[xshift=6.2cm]
  \node[vtx,label=below:{$v$}] (v2) at (0,0) {};
  \node[vtx,label=left:{$u$}] (u2) at (-1.8,0.9) {};
  \node[vtx,label=right:{$w$}] (w2) at (1.8,0.9) {};
  \draw[arr] (u2) -- node[above,font=\small] {$x$} node[below,font=\small] {$g$} (w2);
  \draw[stub,gray] (v2) -- (-1.0,-1.1); \draw[stub,gray] (v2) -- (1.0,-1.1);
  \end{scope}
\end{tikzpicture}
\caption{A balanced pair at $v$: the flow enters $v$ along $e$, from $u$, and leaves $v$ along $e'$, to $w$, with the same value
$x$. The two edges are replaced by one edge $g=uw$ with the value $x$; the other edges at $v$ (gray) are unchanged. The
flows with this balanced pair correspond to the flows of the smaller graph (Lemma~\ref{lem:split}).}
\label{fig:split}
\end{figure}

Statement \textup{(LOCAL)} is needed only for loopless graphs in which every degree is $0$ or at least $4$: loops can be
deleted (Lemma~\ref{lem:loop}), and at a vertex of degree $2$ the two edges form a balanced pair in every flow, by the flow condition.
Section~\ref{sec:local} proves it when every degree is $0$ or $4$; this gives Theorem~\ref{thm:deg4}
(all degrees at most~$4$). In Section~\ref{sec:proofs} a vertex of larger degree is replaced by a triangle, with two of
the edges at the vertex going to each of two vertices of the triangle and the others to the third (Figure~\ref{fig:triangle}). This is
the reduction of \cite[Observation~4.6]{CLSWX}, which shows for $|A|\ge 5$ that all even graphs are $A$-flow-connected
if and only if all $4$-regular graphs are. We need only one direction, and it holds for $|A|\ge 4$: if the graph with
the triangle is $A$-flow-connected, then so is the original graph (Lemma~\ref{lem:triangle}, a special case of
\cite[Theorem~5.9(1)]{CLSWX}, see Remark~\ref{rem:Aconnected}; we include its short proof).

Part~(b) of Theorem~\ref{thm:main} does not follow from part~(a) by a general principle: whether every
$\Z_k$-flow-connected graph is $k$-flow-connected is, for general $k$, an open problem \cite[Problem~2.6]{EHLMNS}. The steps of Section~\ref{sec:local} can have any value, but
the lemmas that turn a step between two nowhere-zero $\Z_k$-flows into a step between two $k$-flows with these residues modulo~$k$
(Lemmas~\ref{lem:step1} and~\ref{lem:step2}) apply only to steps with values $\pm1$
and $\pm2$. For the cyclic groups $\Z_m$ with $m\ge 4$,
Section~\ref{sec:cyclic} therefore proves \textup{(LOCAL)} using only steps with values $\pm1$ and $\pm2$, by a counting argument
that needs no bound on the degrees. This gives Theorem~\ref{thm:cyclic}, which also proves part~(a) for cyclic groups
(Remark~\ref{rem:noncyclic} shows where the counting argument fails for other groups). Section~\ref{sec:int} lifts these chains to
integer flows, using Tutte's lifting theorem \cite{Tutte} in the form of \cite[Theorem~2.3]{EHLMNS}: every nowhere-zero
$\Z_k$-flow is the residue modulo $k$ of a $k$-flow.

\medskip
An expanded version of this paper, with further results, is posted on Zenodo \cite{G}. It bounds the diameter of
$\cF(G,A)$ by a polynomial in the number of edges when $G$ is even and $|A|\ge 4$, and it treats graphs with few
vertices of odd degree. None of its results is used here.

\section{Definitions and results}\label{sec:words}

\paragraph{Graphs.} All graphs are finite; loops and parallel edges are allowed. A graph is \emph{even} if every
vertex has even degree, where a loop counts twice. A \emph{cycle} is a connected subgraph in which every vertex has
degree~$2$; a loop is a cycle with one edge, and two parallel edges form a cycle with two edges. A \emph{directed
cycle} is a cycle together with one of its two directions of traversal. An \emph{Eulerian orientation} of a graph
directs every edge so that at every vertex the in-degree equals the out-degree; the directed edges are its
\emph{arcs}. Every even graph has an Eulerian orientation.

\paragraph{Flows.} Let $A$ be a finite abelian group, and let $G$ be a graph. Fix a
reference direction of every edge. An \emph{$A$-flow} of $G$ is a map $f\colon E(G)\to A$ such that at every vertex
the sum of $f$ over the edges that leave the vertex equals the sum over the edges that enter it. The value of $f$ on an
edge $e$ is $f(e)$ in the reference direction of $e$ and $-f(e)$ in the other direction. From now on, a \emph{flow}
is a nowhere-zero $A$-flow, that is, $f(e)\ne 0$ for every edge $e$.

For a directed cycle $Q$, let $\chi_Q(e)$ be $1$ if $Q$ traverses $e$ in its reference direction, $-1$ if
$Q$ traverses $e$ against it, and $0$ if $e$ is not on $Q$. A \emph{step} is a pair of flows $f$ and
$f'=f+t\chi_Q$, where $Q$ is a directed cycle and $t\ne 0$ is an element of $A$, the \emph{value} of the step. The value is determined only up to sign, as reversing $Q$ or exchanging $f$ and $f'$ replaces $t$ by $-t$. The
graph $\cF(G,A)$ has the flows of $G$ as its vertices, and two flows are adjacent if they form a step. In \cite{CLSWX} two different flows are adjacent if they agree outside a cycle $Q$, and in \cite{EHLMNS} if the support of their difference is a cycle $Q$; both give the same graph $\cF(G,A)$, because the
difference of two such flows is an $A$-flow that is zero outside $Q$, and such an $A$-flow is $t\chi_Q$ for some $t\in A$ and a direction of $Q$. The graph $G$ is \emph{$A$-flow-connected} if $\cF(G,A)$ is connected; here $\cF(G,A)$ counts as connected when it has at most one
vertex, in particular when $G$ has no~flow.

Let $T$ be a set of non-zero elements of $A$ with $T=-T$. A \emph{$T$-step} is a step whose value is in $T$, and two flows
are \emph{$T$-joined} if a chain of $T$-steps (possibly empty) leads from one to the other.

A \emph{multi-step} with value $t\ne 0$ is a pair of flows $f$ and $f'=f+t(\chi_{Q_1}+\dots+\chi_{Q_r})$, where
$r\ge 1$ and $Q_1,\dots,Q_r$ are directed cycles, no two of which have a common edge. Every partial sum
$f+t(\chi_{Q_1}+\dots+\chi_{Q_j})$ is then a flow, as it agrees with $f$ or $f'$ on every edge; so a
multi-step is a chain of $r$ steps with value $t$.

\paragraph{Ends of edges.} Every edge has two ends, one at each endpoint; a loop at $v$ has both ends at $v$. The degree of $v$ is the
number of ends at $v$. If $v$ is a vertex of a cycle $Q$, the \emph{ends of $Q$ at $v$} are the two ends at $v$ of the
edges of $Q$. For an $A$-flow $f$ and an end $i$ at $v$ of the edge $e$, the \emph{outward value}
$y_i(f)\in A$ is the value of $f$ on $e$ in the direction that leaves $v$ through $i$; it is non-zero if $f$ is a flow. The flow
condition at $v$ says
\begin{equation}\label{eq:kirchhoff}
\sum_{i\text{ an end at }v} y_i(f)=0 .
\end{equation}
If a directed cycle $Q$ enters $v$ through the end $i$ and leaves $v$ through the end $j$, then
\begin{equation}\label{eq:ends}
y_i(f+t\chi_Q)=y_i(f)-t,\qquad y_j(f+t\chi_Q)=y_j(f)+t ,
\end{equation}
and the other ends at $v$ keep their outward values.

\begin{definition}\label{def:balanced}
Let $G$ have no loop. A \emph{balanced pair} of a flow $f$ at a vertex $v$ is a pair of different ends $i,j$ at
$v$ with $y_i(f)+y_j(f)=0$; that is, $f$ enters $v$ through one of them and leaves $v$ through the other with the same
value. (Loops are excluded because the two ends of a loop would always form a balanced pair.)
\end{definition}

We prove the following three theorems.
\begin{theorem}\label{thm:deg4}
Let $A$ be a finite abelian group. Every even graph in which all degrees are at most $4$ is $A$-flow-connected.
\end{theorem}

\begin{theorem}\label{thm:all}
Let $A$ be a finite abelian group with $|A|\ge 4$. Every even graph is $A$-flow-connected.
\end{theorem}

\begin{theorem}\label{thm:cyclic}
Let $m\ge 4$. In every even graph, any two nowhere-zero $\Z_m$-flows are joined by a chain of steps whose values are
in $\{1,-1,2,-2\}$.
\end{theorem}

Sections~\ref{sec:proofs} and~\ref{sec:cyclic} prove these theorems, and Section~\ref{sec:int} deduces Theorem~\ref{thm:main}~(b) from Theorem~\ref{thm:cyclic}.

\section{Reduction to flows with a balanced pair}\label{sec:red}

This section proves Proposition~\ref{prop:red}, the reduction described in the outline. In this section $A$ is a finite abelian group with at least two elements, and $T$ is a set of non-zero elements of $A$
with $T=-T$. Here $T$ is the set of allowed step values ($A\setminus\{0\}$ for Theorems~\ref{thm:deg4} and~\ref{thm:all}, and
$\{1,-1,2,-2\}$ for Theorem~\ref{thm:cyclic}). We use two conditions on $T$:
\begin{itemize}
\item[(T1)] any two non-zero elements of $A$ are joined by a sequence of non-zero elements of $A$ in which
consecutive elements differ by an element of $T$;
\item[(T2)] there is a non-zero $a\in A$ with $2a=0$ or $2a\in T$.
\end{itemize}
Both conditions hold for $T=A\setminus\{0\}$, and also for $A=\Z_m$ with $m\ge 3$ and $T=\{1,-1,2,-2\}$: the sequence
$1,2,\dots,m-1$ gives (T1), and $a=1$ gives (T2). Condition (T1) is used for loops in Lemma~\ref{lem:loop}, and (T2) is
used in Proposition~\ref{prop:red} to join the constant flows of Lemma~\ref{lem:const}.

\begin{lemma}[deleting a loop]\label{lem:loop}
Let $T$ satisfy \textup{(T1)}, and let $\ell$ be a loop of $G$. If any two flows of $G-\ell$ are $T$-joined, then any two
flows of $G$ are $T$-joined.
\end{lemma}
\begin{proof}
A loop adds the same term to both sides of the flow condition at its vertex. So a flow of $G$ is a flow of $G-\ell$
together with any non-zero value on $\ell$. A step of $G-\ell$, with the value on $\ell$ kept, is a step of $G$, and a
change of the value on $\ell$ by an element of $T$ is a $T$-step along the cycle $\ell$. Given two flows of $G$, join their
restrictions by $T$-steps of $G-\ell$, then change the value on $\ell$ along a sequence as in \textup{(T1)}.
\end{proof}

\paragraph{Splitting off.} Figure~\ref{fig:split} shows the following operation. Let $G$ have no loop, let $v$ be a
vertex, and let $e=vu$ and $e'=vw$ be two different edges at $v$ ($u=w$ is
allowed). Let $G'$ be the graph $G-e-e'+g$ with a new edge $g=uw$ (a loop if $u=w$): the first end of $g$ is at $u$ and
replaces the end of $e$ at $u$, and its second end is at $w$ and replaces the end of $e'$ at $w$. The direction of
$g$ \emph{from $u$ to $w$} is the direction from its first end to its second end; if $u=w$, this fixes one of the two
directions of the loop $g$. For a map $f'\colon E(G')\to A$, let $x$ be the value of $f'$ on $g$
from $u$ to $w$, and let $Lf'\colon E(G)\to A$ be the map that agrees with $f'$ on the edges other than $e,e',g$ (the common edges of $G$ and $G'$) and has the value $x$ on
$e$ in the direction from $u$ to $v$ and on $e'$ in the direction from $v$ to $w$. Let $B(v;e,e')$ be the set of flows
$f$ of $G$ in which the ends of $e$ and $e'$ at $v$ are a balanced pair.

\begin{lemma}[splitting off]\label{lem:split}
Let $G$, $v$, $e$, $e'$, $G'$, $L$ and $B(v;e,e')$ be as above.
\begin{itemize}
\item[(a)] $L$ is a bijection from the flows of $G'$ onto $B(v;e,e')$.
\item[(b)] If $f'_2$ comes from $f'_1$ by a multi-step of $G'$ with value $t$, then $Lf'_2$ comes from $Lf'_1$ by
a multi-step of $G$ with value $t$. In particular, if two flows of $G'$ are $T$-joined, their images under $L$ are
$T$-joined.
\end{itemize}
\end{lemma}
\begin{proof}
(a) In $Lf'$ the ends of $e$ and $e'$ at $v$ have the outward values $-x$ and $x$, so the flow condition holds at $v$;
at $u$ the edge $e$ leaves with value $x$ as $g$ did, and at $w$ the edge $e'$ enters with value $x$ as $g$
did. So $Lf'$ is a flow in $B(v;e,e')$. Conversely, for $f\in B(v;e,e')$ the value of $f$ on $e$ from $u$ to $v$
equals its value on $e'$ from $v$ to $w$; giving $g$ this value from $u$ to $w$ and keeping the other values defines
the only map $f'$ with $Lf'=f$; it is a flow, as the ends of $g$ have the outward values of the ends of $e,e'$ that they
replace, and the ends of $e,e'$ at $v$ have outward values with sum~$0$.

(b) Let $Q'$ be a directed cycle of $G'$. If $g$ is not on $Q'$, let $W=Q'$, a directed cycle of $G$. If $g$ is on
$Q'$, replace $g$ by the two edges $e,e'$, traversed through $v$ in the direction of $Q'$ (first $e$ from $u$ to $v$
and then $e'$ from $v$ to $w$ if $Q'$ traverses $g$ from $u$ to $w$, and the reverse otherwise); this gives a closed
walk $W$ of $G$ that uses every edge at most once. If $u=w$, then $Q'=g$ and $W$ is the cycle formed by $e$ and $e'$.
If $u\ne w$ and $Q'$ does not pass through $v$, then $W$ is a directed cycle. If $u\ne w$ and $Q'$ passes through $v$,
then $W$ visits $v$ twice and no other vertex twice, and splitting $W$ at $v$ gives two directed cycles with no
common edge. In every case let $\chi_W$ be the sum of the maps $\chi$ of these one or two directed cycles; then
$L(f'+t\chi_{Q'})=Lf'+t\chi_W$. If $Q'_1,\dots,Q'_r$ are directed cycles of $G'$, no two of which have a common edge,
then no two of the directed cycles of $G$ obtained from them have a common edge: they use the edges of the $Q'_j$,
with $e,e'$ instead of $g$. So $L$ maps a multi-step of $G'$ to a multi-step of $G$ with the same value.
\end{proof}

Part~(a) of the next lemma restates \cite[Lemma~2.3]{CLSWX}, which is about Eulerian orientations, in the language of
flows; we include its short proof.

\begin{lemma}[constant flows]\label{lem:const}
Let $G$ be even and let $a\ne 0$ be an element of $A$. For an Eulerian orientation $D$ of $G$, let
$a\one_D$ be the constant flow with value $a$ along every arc of $D$.
\begin{itemize}
\item[(a)] For two Eulerian orientations $D_1,D_2$ of $G$, either $a\one_{D_2}=a\one_{D_1}$, or $a\one_{D_2}$ comes from
$a\one_{D_1}$ by a multi-step with value $-2a$.
\item[(b)] In the situation of Lemma~\ref{lem:split}, $B(v;e,e')$ contains $a\one_D$ for some Eulerian orientation $D$ of $G$.
\end{itemize}
\end{lemma}
\begin{proof}
(a) The arcs of $D_1$ that are reversed in $D_2$ form a directed graph with in-degree equal to out-degree at every
vertex, because $D_1$ and $D_2$ are both Eulerian orientations; so they are the union of directed cycles $Q_1,\dots,Q_r$ of $D_1$, no two of which have a common arc. The
reversal of an arc changes its value along $D_1$ from $a$ to $-a$. So
$a\one_{D_2}=a\one_{D_1}-2a(\chi_{Q_1}+\dots+\chi_{Q_r})$. If $2a=0$ or $r=0$, the two flows are equal; otherwise this is a
multi-step with value $-2a$.

(b) The graph $G'$ is even. Take an Eulerian orientation $D'$ of $G'$, and let $D$ agree with $D'$ on the common edges and
direct $e$ and $e'$ through $v$ in the direction of $g$ in $D'$ (from $u$ to $v$ to $w$ if $D'$ directs $g$ from $u$
to $w$, and the reverse otherwise). Then $D$ is an Eulerian orientation of $G$ and $a\one_D=L(a\one_{D'})$.
\end{proof}

\begin{proposition}[reduction to a balanced pair]\label{prop:red}
Let $\cC$ be a class of even graphs such that $G-\ell\in\cC$ for every $G\in\cC$ and every loop $\ell$ of $G$, and
$G'\in\cC$ for every $G\in\cC$ with no loop and every graph $G'$ obtained from $G$ as in Lemma~\ref{lem:split}. Let
$T$ satisfy \textup{(T1)} and \textup{(T2)}. Suppose that the following statement holds.
\begin{itemize}[leftmargin=4.6em,labelsep=0.5em]
\item[\textup{(LOCAL)}] For every loopless $G\in\cC$ in which every vertex has degree $0$
or at least $4$, and which has at least one edge, every flow of $G$ is $T$-joined to a flow that has a balanced pair at some vertex.
\end{itemize}
Then in every $G\in\cC$ any two flows are $T$-joined.
\end{proposition}
\begin{proof}
We use induction on the number of edges of $G$. A graph with no edge has exactly one flow. If $G$ has a loop $\ell$, use
Lemma~\ref{lem:loop} and the induction hypothesis for $G-\ell$. From now on $G$ has no loop.

If some vertex $v$ has degree $2$, its two ends belong to two different edges $e,e'$, and by \eqref{eq:kirchhoff}
they are a balanced pair in every flow. So every flow of $G$ is in $B(v;e,e')$, and by Lemma~\ref{lem:split} and the
induction hypothesis for $G'$, any two flows of $G$ are $T$-joined.

Otherwise $G$ has an edge and every vertex has degree $0$ or at least $4$. Fix $a$ as in \textup{(T2)}. For every
vertex $v$ and any two different edges $e,e'$ at $v$, any two flows in $B(v;e,e')$ are $T$-joined, by the induction
hypothesis for $G'$ and Lemma~\ref{lem:split}, and $B(v;e,e')$ contains a flow $a\one_D$, by
Lemma~\ref{lem:const}~(b). By Lemma~\ref{lem:const}~(a) and \textup{(T2)}, all flows $a\one_D$ are $T$-joined to one another (if $2a\ne 0$, the value
$-2a$ is in $T$, because $2a\in T$ and $T=-T$). A flow with a balanced pair is in some $B(v;e,e')$, because $G$ has
no loop. So all flows with a balanced pair are $T$-joined to one another, and \textup{(LOCAL)} joins every flow to
one of them.
\end{proof}

\section{Reaching a balanced pair when all degrees are 0 or 4}\label{sec:local}

This section proves statement \textup{(LOCAL)} of Proposition~\ref{prop:red} for graphs in which every
degree is $0$ or $4$ (Proposition~\ref{prop:local4}). We fix generators of the cyclic summands of $A$ and call an edge
small if its value is one of these generators or its negative (precise definitions follow). First, steps along cycles
change a flow into one in which the edges that are not small form a forest (Corollary~\ref{cor:forest}); then, as the graph has an edge and every vertex with an edge has degree $4$, the small
edges contain a cycle, and one more step along that cycle gives a balanced pair (Lemma~\ref{lem:Bstar}).

Let $A$ be a finite abelian
group with at least two elements.
Write $A$ as a direct sum of cyclic groups,
\[
A=\langle s_1\rangle\oplus\dots\oplus\langle s_q\rangle,\qquad s_h\ne 0 ,
\]
so that every element $c$ of $A$ can be written in exactly one way as $c=c_1+\dots+c_q$ with $c_h\in\langle
s_h\rangle$; the $c_h$ are the \emph{parts} of $c$. Let $S=\{s_1,-s_1,\dots,s_q,-s_q\}$. The elements of $S$ are
called \emph{small}, and the \emph{color} of $s\in S$ is the index $h$ with $s\in\langle s_h\rangle$. We have $S=-S$ and $0\notin S$, and $S$ generates $A$. For a
flow $f$, an edge is \emph{small} if its value is small (this does not depend on the direction); the \emph{color} of
a small edge is the color of its value, which is the same in both directions. Let $\rho(f)$ be the number of edges
that are not small.

\begin{lemma}[a cycle with no small edge]\label{lem:Astar}
Let $f$ be a flow of a graph $G$ and $Q$ a directed cycle of $G$ none of whose edges is small. Then there is
a step $f'=f+t\chi_Q$ with $\rho(f')\le\rho(f)-1$.
\end{lemma}
\begin{proof}
Let $X$ be the set of values of the edges of $Q$ in the direction of $Q$. Then $0\notin X$
and $X\cap S=\emptyset$. We claim that there are $x\in X$ and $s\in S$ with $x-s\notin X$. If not, then $X-S\subseteq
X$, that is, $X+S\subseteq X$ as $S=-S$; adding elements of $S$ repeatedly gives $X+\langle S\rangle\subseteq X$, and as $S$ generates $A$ and $X$ is not empty,
$X=A$, contradicting $0\notin X$. Let $t=s-x$. Then $t\ne 0$ because $x\notin S$. A new value $x'+t$ with $x'\in X$ is $0$
only if $x'=x-s$, which is not in $X$. So $f'=f+t\chi_Q$ is a flow. An edge of $Q$ with value $x$ has the new
value $s$ and is small in $f'$; the edges outside $Q$ did not change; and no edge of $Q$ was small in $f$. So $\rho(f')\le\rho(f)-1$.
\end{proof}

\begin{corollary}\label{cor:forest}
Every flow $f$ of a graph $G$ is joined by a chain of at most $\rho(f)$ steps to a flow in which the edges that are
not small form a forest.
\end{corollary}
\begin{proof}
Apply Lemma~\ref{lem:Astar} as long as the edges that are not small contain a cycle. Each application is a single step and lowers $\rho$ by at least $1$, and $\rho\ge 0$.
\end{proof}

\begin{lemma}[a cycle of small edges]\label{lem:Bstar}
Let $G$ be a loopless graph in which every vertex has degree $0$ or $4$. Let $f$ be a flow of $G$ with no balanced
pair, and $Q$ a cycle of $G$ all of whose edges are small. Then there are a direction of $Q$ and an element $t\ne 0$ of $A$ such that $f+t\chi_Q$ is a flow with a balanced pair.
\end{lemma}
\begin{proof}
Let $v$ be a vertex of $Q$, to be chosen below. As $v$ lies on $Q$, its degree is $4$; number the four ends at $v$ as $1,2,3,4$, where $1$ and $2$ are the ends at $v$ of the two edges of $Q$ at $v$, and write
$y_i=y_i(f)$. Direct $Q$ so that it enters $v$ through the end $1$ and leaves through the end $2$. For $j\in\{3,4\}$
let $t_j=y_1+y_j$, and let $j'$ be the other element of $\{3,4\}$; $t_j$ is not $0$ because $f$ has no balanced pair. By \eqref{eq:ends} and
\eqref{eq:kirchhoff}, in $f+t_j\chi_Q$ the end $1$ has the outward value $y_1-t_j=-y_j$ and the end $2$ has the
outward value $y_2+t_j=-y_{j'}$. Neither is $0$, and the ends $1$ and $j$ are a balanced pair. An edge of $Q$ that is not at $v$ gets
the value $0$ exactly if its value in the direction of $Q$ was $-t_j$. Hence the following holds.
\begin{itemize}
\item[$(\ast)$] If no edge of $Q$ that is not at $v$ has the value $-t_j$ in the direction of $Q$, then
$f+t_j\chi_Q$ is a flow with a balanced pair.
\end{itemize}
In particular, $(\ast)$ applies if $t_j$ is not small, because then $-t_j$ is not small either ($S=-S$), while all
values on $Q$ are small. Below
we use the identity $t_3+t_4=2y_1+y_3+y_4=y_1-y_2$, which holds by \eqref{eq:kirchhoff}.

\emph{Case 1: two consecutive edges of $Q$ have different colors.} Choose $v$ to be a common vertex of these two edges. Then $y_1=a$ is
small, of color $h$, and $y_2=b$ is small, of color $h'\ne h$. Suppose that $t_3$ and $t_4$ are both small. Their sum
$t_3+t_4=a-b$ has exactly two non-zero parts, $a$ in $\langle s_h\rangle$ and $-b$ in $\langle s_{h'}\rangle$. If $t_3$
and $t_4$ had the same color, both would lie in one cyclic summand of $A$, and their sum would have at most one non-zero
part; so they have different colors. Then $t_3$ and $t_4$ are non-zero elements of two different cyclic summands of $A$,
so by the uniqueness of the parts they are the two non-zero parts of $a-b$, that is, $\{t_3,t_4\}=\{a,-b\}$. So $t_3=a$ or $t_4=a$; but $t_3=a$ means $y_3=0$, and $t_4=a$ means $y_4=0$. This is
impossible. So $t_3$ or $t_4$ is not small, and $(\ast)$ gives the step.

\emph{Case 2: all edges of $Q$ have the same color $h$.} Here $v$ may be any vertex of $Q$. Then $y_1,y_2\in\{s_h,-s_h\}$
and $y_1+y_2\ne 0$, so $y_1=y_2=a$ and $2a\ne 0$. (If $2s_h=0$ this case does not occur: the ends $1,2$ would be a
balanced pair.) The values of the edges of $Q$ in the direction of $Q$ are $a$ or $-a$. Now $-t_3=a$ means
$y_3=-2a$ and then $y_4=-2a-y_3=0$; and $-t_3=-a$ means $y_3=0$. Both are impossible, so $(\ast)$ applies with $j=3$ and gives a step with value $t_3$.
\end{proof}

\begin{proposition}[reaching a balanced pair when all degrees are $0$ or $4$]\label{prop:local4}
Let $G$ be a loopless graph in which every vertex has degree $0$ or $4$, and which has at least one edge. Every flow $f$ of $G$ is joined by a chain of at most $\rho(f)+1$ steps to a flow with a balanced pair.
\end{proposition}
\begin{proof}
By Corollary~\ref{cor:forest}, $f$ is joined by a chain of at most $\rho(f)$ steps to a flow $f^*$ in which the
edges that are not small form a forest. Let $n\ge 1$ be the number of vertices of degree $4$. Then $G$ has $2n$ edges, all on these
$n$ vertices; at most $n-1$ of them are not small in $f^*$ (a forest on $n$ vertices has at most $n-1$ edges), so at
least $n+1$ are small, and the small edges contain a cycle $Q$. If $f^*$ has no balanced pair,
Lemma~\ref{lem:Bstar} gives a further step to a flow with a balanced pair.
\end{proof}

\section{Finite abelian groups: proofs of Theorems \ref{thm:deg4} and \ref{thm:all}}\label{sec:proofs}

\begin{proof}[Proof of Theorem~\ref{thm:deg4}]
If $|A|=1$, then $G$ has at most one flow (none if it has an edge), so $G$ is $A$-flow-connected. Let $|A|\ge 2$, let $\cC$ be the class of even graphs with all degrees at most $4$, and let $T=A\setminus\{0\}$; it satisfies \textup{(T1)} and \textup{(T2)}. Deleting a loop lowers the degree of its vertex by $2$. Splitting off as in Lemma~\ref{lem:split} lowers the degree of $v$ by $2$ and keeps the other degrees
(a new loop at $u=w$ counts twice). So $\cC$ is closed under both operations. A graph in $\cC$ with no loop and no
vertex of degree $2$ has only vertices of degree $0$ or $4$, so \textup{(LOCAL)} is Proposition~\ref{prop:local4}.
Proposition~\ref{prop:red} gives the theorem.
\end{proof}

For the proof of Theorem~\ref{thm:all}, a vertex of degree at least $6$ is replaced by a triangle
(Figure~\ref{fig:triangle}); Lemma~\ref{lem:triangle} shows that if the new graph is
$A$-flow-connected, then so is the original graph.

\begin{figure}[ht]
\centering
\begin{tikzpicture}[scale=1.0]
  \node[vtx] (v) at (0,0) {};
  \node[font=\small,fill=white,inner sep=1pt] at (0.34,0.14) {$v$};
  \foreach \ang in {0,45,...,315} \draw[stub] (v) -- (\ang:1.3);
  \node[font=\small] at (0,-1.9) {degree $8$};
  \draw[-{Stealth[length=3mm]},very thick] (2.0,0) -- (3.0,0);
  \begin{scope}[xshift=5.8cm]
  \node[vtx] (x1) at (-1.0,-0.55) {};
  \node[vtx] (x2) at (1.0,-0.55) {};
  \node[vtx] (x3) at (0,0.9) {};
  \node[font=\small] at (-1.26,-0.26) {$x_1$};
  \node[font=\small] at (1.26,-0.26) {$x_2$};
  \node[font=\small] at (0,0.36) {$x_3$};
  \draw[thick] (x1) -- (x2) -- (x3) -- (x1);
  \draw[stub] (x1) -- ++(190:1.1); \draw[stub] (x1) -- ++(250:1.1);
  \draw[stub] (x2) -- ++(-10:1.1); \draw[stub] (x2) -- ++(-70:1.1);
  \foreach \ang in {30,70,110,150} \draw[stub] (x3) -- ++(\ang:1.1);
  \node[font=\small] at (0,-1.9) {degrees $4$, $4$, $6$};
  \end{scope}
\end{tikzpicture}
\caption{A vertex $v$ of degree $d\ge 6$ (here $d=8$) is replaced by a triangle $x_1x_2x_3$: two of the edges at $v$
go to $x_1$, two go to $x_2$, and the other $d-4$ go to $x_3$. The new degrees are $4$, $4$ and $d-2$. This is the
replacement used in the proof of \cite[Observation~4.6]{CLSWX}.}
\label{fig:triangle}
\end{figure}

\begin{lemma}[contracting a triangle]\label{lem:triangle}
Let $|A|\ge 4$. Let $G^+$ be a graph and $\Delta$ a triangle of $G^+$ with vertices $x_1,x_2,x_3$, such that the
three edges of $\Delta$ are the only edges of $G^+$ whose two ends are both at vertices of $\Delta$. Let $G=G^+/\Delta$ be the
graph obtained by deleting the three edges of $\Delta$ and identifying $x_1,x_2,x_3$ to one vertex $x$. For a flow
$f$ of $G^+$, let $pf$ be its restriction to $E(G)=E(G^+)\setminus E(\Delta)$.
\begin{itemize}
\item[(a)] $pf$ is a flow of $G$, and every flow of $G$ is $pf$ for some flow $f$ of $G^+$.
\item[(b)] If $f_2$ comes from $f_1$ by a multi-step with value $t$, then $pf_2=pf_1$ or $pf_2$ comes from $pf_1$
by a multi-step with value $t$.
\end{itemize}
So if $G^+$ is $A$-flow-connected, then $G$ is $A$-flow-connected.
\end{lemma}
\begin{proof}
(a) The flow condition of $pf$ at $x$ is the sum of the flow conditions of $f$ at $x_1,x_2,x_3$, in which the edges
of $\Delta$ cancel. Let $\varphi$ be a flow of $G$. For $i=1,2,3$, let $\beta_i$ be the sum of the outward values of $\varphi$
over the ends at $x$ that are ends at $x_i$ in $G^+$; then $\beta_1+\beta_2+\beta_3=0$. Give the edges $x_1x_2$,
$x_2x_3$, $x_3x_1$ of $\Delta$, in these directions, the values $z$, $z-\beta_2$, $z+\beta_1$ (the value $z$ on
$x_1x_2$ is free, and the flow conditions at $x_2$ and $x_1$ then fix the other two values). The flow conditions at
$x_1,x_2,x_3$ are $\beta_1+z-(z+\beta_1)=0$, $\beta_2+(z-\beta_2)-z=0$ and $\beta_3+(z+\beta_1)-(z-\beta_2)=0$; they
hold for every $z$. The three values are not $0$ if $z\notin\{0,\beta_2,-\beta_1\}$, and $|A|\ge 4$ leaves such a
$z$.

(b) Let $f_2=f_1+t(\chi_{Q_1}+\dots+\chi_{Q_r})$ with directed cycles of $G^+$, no two of which have a common edge.
Let $H$ be the directed graph on $V(G)$ whose arcs are the edges of $Q_1,\dots,Q_r$ outside $\Delta$, each with the
direction of its cycle. At every vertex other than $x$, the in-degree of $H$ equals its out-degree, because this holds for the union of the
cycles. It also holds at $x$: in the union of the cycles, the total in-degree of $x_1,x_2,x_3$ equals their total
out-degree, and every arc of the cycles inside $\Delta$ adds one to the total in-degree and one to the total out-degree. So $H$ is the union of
directed cycles of $G$, no two of which have a common arc, and $pf_2-pf_1$ is $t$ times the sum of their maps $\chi$.
If $H$ has no arc, $pf_2=pf_1$; otherwise this is a multi-step of $G$ with value $t$.

The last statement follows, because a multi-step is a chain of steps and $p$ is onto by (a).
\end{proof}

\begin{remark}\label{rem:Aconnected}
Lemma~\ref{lem:triangle} is the special case of \cite[Theorem~5.9(1)]{CLSWX} for a triangle, which is $A$-connected
in the sense of \cite{JLPT} when $|A|\ge 4$ (this notion is different from $A$-flow-connectedness).
\end{remark}

\begin{proof}[Proof of Theorem~\ref{thm:all}]
By Lemma~\ref{lem:loop} (with $T=A\setminus\{0\}$), deleting the loops one by one, we may assume that $G$ has no loop. Let
$\Phi(G)=\sum_v\max(\deg v-4,0)$. We use induction on $\Phi(G)$. If $\Phi(G)=0$, this is Theorem~\ref{thm:deg4}.
Otherwise let $v$ be a vertex of degree $d\ge 6$. Let $G^+$ be the graph obtained from $G$ by replacing $v$ with a
triangle $\Delta$ with vertices $x_1,x_2,x_3$: two of the ends at $v$ become ends at $x_1$, two become ends at $x_2$,
and the other $d-4$ become ends at $x_3$ (Figure~\ref{fig:triangle}). The new degrees are $4$, $4$ and $d-2$, so $G^+$ is even and
$\Phi(G^+)=\Phi(G)-2$. The graph $G^+$ has no loop, and the edges of $\Delta$ are the only edges whose two ends
are both at vertices of $\Delta$, because $G$ has no loop at $v$. Also $G=G^+/\Delta$. By the induction hypothesis $G^+$ is
$A$-flow-connected, and Lemma~\ref{lem:triangle} gives the same for $G$.
\end{proof}

\section{Cyclic groups: steps with values \texorpdfstring{$\pm1$}{+-1} and \texorpdfstring{$\pm2$}{+-2}}\label{sec:cyclic}

This section proves Theorem~\ref{thm:cyclic}. Let $A=\Z_m$ with $m\ge 4$ and $T=\{1,-1,2,-2\}$. An edge is small for a flow $f$ if its value is $1$
or $-1$ (as in Section~\ref{sec:local} with $S=\{1,-1\}$), and $\rho(f)$ is the number of edges that are not small. The value of an edge that is not
small is, in each direction, the residue of an integer in $\{2,\dots,m-2\}$.

\begin{lemma}[a cycle with no small edge, in $\Z_m$]\label{lem:Acyc}
Let $f$ be a flow of a graph $G$ and $Q$ a cycle none of whose edges is small. Then $f$ is $T$-joined, by steps along
$Q$ with value $-1$, to a flow $f'$ with $\rho(f')\le\rho(f)-1$.
\end{lemma}
\begin{proof}
Fix a direction of $Q$. For an edge $e$ of $Q$, let $x_e\in\{2,\dots,m-2\}$ be the integer whose residue is the value
of $e$ in this direction, and let $x$ be the least of the $x_e$. For $j=1,\dots,x-1$ the map $f-j\chi_Q$ has on $e$
the residue of $x_e-j$, where $1\le x_e-j\le m-3$; so it is a flow. In $f'=f-(x-1)\chi_Q$, an edge with $x_e=x$
has the value $1$; the edges outside $Q$ did not change, and no edge of $Q$ was small in $f$. So $\rho(f')\le\rho(f)-1$.
\end{proof}

\begin{lemma}[a cycle of small edges, in $\Z_m$]\label{lem:Bcyc}
Let $G$ have no loop, let $f$ be a flow with no balanced pair, $Q$ a cycle all of whose edges are small, and $v$ a
vertex of $Q$.
\begin{itemize}
\item[(a)] The two ends of $Q$ at $v$ have the same outward value $a\in\{1,-1\}$.
\item[(b)] If an end at $v$ that is not an end of $Q$ has an outward value $c\ne -2a$, then $f$ is $T$-joined, by steps along $Q$, to a
flow with a balanced pair at $v$.
\end{itemize}
\end{lemma}
\begin{proof}
(a) Both values are in $\{1,-1\}$ and their sum is not $0$.

(b) Call the two ends of $Q$ at $v$ the end $1$ and the end $2$, direct $Q$ so that it enters $v$ through the end
$1$ and leaves through the end $2$, and let $t=a+c$. Then $t\ne
0$ ($f$ has no balanced pair), $t\ne a$ ($c\ne 0$), and $t\ne -a$ ($c\ne -2a$). So $t$ is the residue of an integer
$\tau$ with $2\le\tau\le m-2$. The values of the edges of $Q$ in the direction of $Q$ are $1$ or $-1$. For
$j=2,3,\dots,\tau$ the map $f+j\chi_Q$ has on the edges of $Q$ the residues of $1+j$ or $-1+j$, integers between $1$ and
$m-1$; so it is a flow. The chain $f$, $f+2\chi_Q$, $f+3\chi_Q$, \dots, $f+\tau\chi_Q$ consists of one step with
value $2$ and $\tau-2$ steps with value $1$. (The chain does not start with a step with value $1$: by (a), one of
the two edges of $Q$ at $v$ has the value $-1$ in the direction of $Q$, and $f+\chi_Q$ is $0$ on that edge.) In $f+\tau\chi_Q=f+t\chi_Q$ the end $1$ has the outward value $a-t=-c$ by
\eqref{eq:ends}, and the end with the outward value $c$ did not change; so these two ends are a balanced pair at $v$.
\end{proof}

\begin{proposition}[reaching a balanced pair in $\Z_m$]\label{prop:localcyc}
Let $m\ge 4$, and let $G$ be a loopless even graph in which every vertex has degree $0$ or at least $4$, and which has at least one edge. Every flow of $G$ is $T$-joined to a flow with a balanced pair.
\end{proposition}
\begin{proof}
Let $n\ge 1$ be the number of vertices of degree at least $4$; all edges are on these vertices, and $|E(G)|\ge 2n$.
By Lemma~\ref{lem:Acyc}, used as long as the edges that are not small contain a cycle, a given flow is $T$-joined to a flow
$f^*$ in which they form a forest; so at most $n-1$ edges are not small, as they lie on the $n$ vertices. Suppose that
$f^*$ has no balanced pair (otherwise there is nothing to prove), and that the hypothesis of Lemma~\ref{lem:Bcyc}~(b) fails for every cycle $Q$ of small
edges and every vertex $v$ of $Q$. Then at every vertex $v$ of such a cycle, all ends at $v$ that are not ends of the cycle have the outward
value $-2a$, where $a\in\{1,-1\}$ is the common outward value at $v$ of the two ends of the cycle
(Lemma~\ref{lem:Bcyc}~(a)). As $m\ge 4$, $-2a$ is neither $1$ nor $-1$: $-2a=a$ would mean $3\equiv 0$, and $-2a=-a$ would
mean $1\equiv 0$, modulo $m$. So the edges of these ends are not small, and every vertex of a cycle of small edges lies
on exactly two small edges, the two edges of the cycle. Hence every cycle of small edges is a whole component of the graph formed by the small edges on the $n$
vertices. So every component of this graph is a tree or a cycle (a component that contains a cycle is that cycle), and
there are at most $n$ small edges. Then
$|E(G)|\le n+(n-1)<2n$, a contradiction. So, for some cycle $Q$ of small edges and some vertex $v$ of $Q$,
Lemma~\ref{lem:Bcyc}~(b) joins $f^*$ to a flow with a balanced pair.
\end{proof}

\begin{proof}[Proof of Theorem~\ref{thm:cyclic}]
Apply Proposition~\ref{prop:red} to the class $\cC$ of all even graphs and $T=\{1,-1,2,-2\}$: \textup{(T1)} and
\textup{(T2)} hold, as noted at the start of Section~\ref{sec:red}, and \textup{(LOCAL)} is
Proposition~\ref{prop:localcyc}.
\end{proof}

\begin{remark}\label{rem:noncyclic}
Theorem~\ref{thm:cyclic} proves Theorem~\ref{thm:all} again for cyclic groups, without Lemma~\ref{lem:triangle}. The
counting argument of Proposition~\ref{prop:localcyc} uses that the value $-2a$ in Lemma~\ref{lem:Bcyc}~(b) is not
small. This holds in $\Z_m$ with $m\ge 4$, but fails,
for example, in $\Z_3\times\Z_3$ with $S$ as in Section~\ref{sec:local}, where $-2a=a$: then a vertex of a cycle of
small edges can lie on more than two small edges, and the bound of at most $n$ small edges in that argument no longer follows. For example, a vertex of degree $6$ can have six ends
with the same small outward value $a$: the flow condition holds, as $6a=0$; there is no balanced pair, as $2a\ne 0$; and all six edges are small.
\end{remark}

\section{Integer flows: proof of Theorem~\ref{thm:main}~(b)}\label{sec:int}

This section turns a chain of nowhere-zero $\Z_k$-flows from Theorem~\ref{thm:cyclic}
into a walk in $\cF(G,k)$: every $\Z_k$-flow of the chain is replaced by a path of $k$-flows whose residues modulo $k$
are that $\Z_k$-flow. Throughout this section, a \emph{$\Z_k$-flow} means a nowhere-zero $\Z_k$-flow (a flow in the
sense of Section~\ref{sec:words} for $A=\Z_k$); we name the group to tell these flows apart from the $k$-flows, whose
values are integers.

\paragraph{Definitions.} An \emph{integer flow} of a graph $G$ is a map $F\colon E(G)\to\Z$ that satisfies the flow condition at every vertex. For a map $x$ from $E(G)$ to $\Z$ or to the real numbers, the \emph{boundary}
$(\partial x)(v)$ at a vertex $v$ is the sum of $x$ over the edges that enter $v$ minus the sum over the edges that
leave $v$ (a loop at $v$ is in both sums). So an integer map $x$ is an integer flow if and only if $\partial x=0$; also $\partial\chi_Q=0$ for every directed cycle~$Q$.

Let $k\ge 2$. The $k$-flows, the graph $\cF(G,k)$ and $k$-flow-connectedness are defined in Section~\ref{sec:intro};
in the notation of Section~\ref{sec:words}, two $k$-flows are adjacent if they differ by $\lambda\chi_Q$, where $Q$ is a
directed cycle and $\lambda\ne 0$ is an integer. As for $\cF(G,A)$, the graph $\cF(G,k)$ counts as connected when it has
at most one vertex, and it is the same as the graph $\cF(G,k)$ of
\cite{CLSWX,EHLMNS}. An \emph{integer lift} of a $\Z_k$-flow $f$ is a $k$-flow $F$ with $F(e)\equiv f(e)$
modulo $k$ for every edge $e$. The residues modulo $k$ of a $k$-flow form a $\Z_k$-flow, because no value of a
$k$-flow is a multiple of $k$.

\begin{theorem}\label{thm:int}
Let $k\ge 4$. Every even graph is $k$-flow-connected.
\end{theorem}

This is part~(b) of Theorem~\ref{thm:main}. The proof lifts a chain from Theorem~\ref{thm:cyclic} step by step.
Lemma~\ref{lem:lift} is Tutte's lifting theorem in the form of \cite[Theorem~2.3]{EHLMNS}, and Lemma~\ref{lem:fibre}
is contained in the proof of \cite[Lemma~5.5]{EHLMNS}. We include their short proofs to keep the paper self-contained; the rounding lemma
(Lemma~\ref{lem:round}) gives Lemma~\ref{lem:lift} and is used again in the proof of Lemma~\ref{lem:step2}.

\begin{lemma}[rounding]\label{lem:round}
Let $G$ be a graph and $x$ a map from $E(G)$ to the real interval $[0,1]$ such that $(\partial x)(v)$ is an integer
for every vertex $v$. Then there is a map $z\colon E(G)\to\{0,1\}$ with $\partial z=\partial x$ and with $z(e)=x(e)$
for every edge $e$ with $x(e)\in\{0,1\}$.
\end{lemma}
\begin{proof}
Call an edge $e$ \emph{fractional} if $0<x(e)<1$. If no edge is fractional, take $z=x$. Otherwise the fractional
edges contain a cycle $Q$. If some loop is fractional, it is such a cycle. If no loop is fractional, then no vertex $v$
meets exactly one fractional edge, because the other edges at $v$ have integer values and $(\partial x)(v)$ is an
integer; and a graph with no loop, at least one edge and no vertex of degree $1$ contains a cycle. Give $Q$ a direction, and
let $\delta>0$ be the least number such that $x+\delta\chi_Q$ has the value $0$ or $1$ on some edge of $Q$. The map
$x+\delta\chi_Q$ takes values in $[0,1]$, has the same boundary as $x$, agrees with $x$ on the edges that are not
fractional, and has fewer fractional edges than $x$. After at most $|E(G)|$ such changes no edge is fractional; take
$z$ to be the final map.
\end{proof}

For a $\Z_k$-flow $f$ and an edge $e$, let $r_f(e)\in\{1,\dots,k-1\}$ be the integer whose residue is the value of
$f$ on $e$ in the reference direction. Then $\partial r_f\equiv 0$ modulo $k$ at every vertex, and an integer lift of
$f$ has on every edge $e$ the value $r_f(e)$ or $r_f(e)-k$.

\begin{lemma}[integer lifts]\label{lem:lift}
Let $k\ge 2$. Every $\Z_k$-flow $f$ of a graph $G$ has an integer lift.
\end{lemma}
\begin{proof}
The map $x=r_f/k$ takes values in $(0,1)$ and has integer boundary. Let $z$ be as in Lemma~\ref{lem:round}, and
$F=r_f-kz$. Then $\partial F=\partial r_f-k\,\partial z=\partial r_f-k\,\partial x=0$, $F\equiv f$ modulo $k$, and
$F(e)$ is $r_f(e)\in[1,k-1]$ or $r_f(e)-k\in[-(k-1),-1]$.
\end{proof}

\begin{lemma}[integer lifts of a given flow]\label{lem:fibre}
Let $k\ge 2$ and let $f$ be a $\Z_k$-flow of a graph $G$. Any two integer lifts of $f$ are joined by a path in
$\cF(G,k)$ all of whose vertices are integer lifts of $f$.
\end{lemma}
\begin{proof}
Let $F,F'$ be integer lifts of $f$. On every edge both have the value $r_f(e)$ or $r_f(e)-k$, so $\eta=(F'-F)/k$ is an
integer flow with values in $\{-1,0,1\}$. Direct every edge $e$ with $\eta(e)\ne 0$ so that $\eta$ is $1$ in this
direction. By the flow condition these arcs form a directed graph with in-degree equal to out-degree at every vertex,
so they are the union of directed cycles $Q_1,\dots,Q_j$, no two of which have a common arc, and
$\eta=\chi_{Q_1}+\dots+\chi_{Q_j}$. Let $F_i=F+k(\chi_{Q_1}+\dots+\chi_{Q_i})$. On every edge $F_i$ has the value of $F$
or the value of $F'$, so $F_i$ is an integer lift of $f$; and $F_i-F_{i-1}=k\chi_{Q_i}$. So $F=F_0,F_1,\dots,F_j=F'$
is a path in $\cF(G,k)$ of integer lifts of $f$.
\end{proof}

\begin{lemma}[a step with value $1$]\label{lem:step1}
Let $k\ge 2$, let $Q$ be a directed cycle of a graph $G$, and let $f$ and $f+\chi_Q$ be $\Z_k$-flows of $G$. Then for
every integer lift $F$ of $f$, the map $F+\chi_Q$ is an integer lift of $f+\chi_Q$.
\end{lemma}
\begin{proof}
Choose the reference directions so that $Q$ traverses all its edges in the reference direction. (The statement does
not depend on this choice: a change of the reference direction of an edge changes the signs of $F(e)$, of $f(e)$ and
of $\chi_Q(e)$.) The map $F+\chi_Q$ is an integer flow with the residues of $f+\chi_Q$. Let $e$ be an edge of $Q$. Then
$r_f(e)\ne k-1$, because $f+\chi_Q$ is not $0$ on $e$. The value $F(e)$ is $r_f(e)$ or $r_f(e)-k$, so $F(e)+1$ is in
$[2,k-1]$ or in $[2-k,-1]$.
\end{proof}

\begin{lemma}[a step with value $2$]\label{lem:step2}
Let $k\ge 3$, let $Q$ be a directed cycle of a graph $G$, and let $f$ and $f+2\chi_Q$ be $\Z_k$-flows of $G$. Then
there is an integer lift $F$ of $f$ such that $F+2\chi_Q$ is an integer lift of $f+2\chi_Q$.
\end{lemma}
\begin{proof}
Lemma~\ref{lem:step1} cannot simply be used twice, because $f+\chi_Q$ can have zeros. As in the proof of Lemma~\ref{lem:step1}, choose the reference directions so that $Q$ traverses all its edges in the reference
direction. No edge of $Q$ has $r_f(e)=k-2$, because $f+2\chi_Q$ is not $0$ there. Let $E_1$ be the set of edges $e$ of $Q$ with $r_f(e)=k-1$, and let
\[
x=(r_f+\chi_Q)/k .
\]
All values of $x$ are in $(0,1]$, and $x(e)=1$ exactly for $e\in E_1$. (The term $\chi_Q$ is there for the edges of
$E_1$: it forces $z(e)=1$ below, so that the lift has the value $-1$ on these edges and the step leads to the value $1$,
not to $k+1$.) As $\partial\chi_Q=0$, the boundary $\partial x=(\partial r_f)/k$ is
an integer at every vertex. By Lemma~\ref{lem:round} there is $z\colon E(G)\to\{0,1\}$ with $\partial z=\partial x$
and with $z(e)=1$ for every $e\in E_1$. Let $F=r_f-kz$. As in the proof of Lemma~\ref{lem:lift}, $F$ is an integer
lift of $f$.

Now $F+2\chi_Q$ is an integer flow with the residues of $f+2\chi_Q$, and outside $Q$ it equals $F$. If $e\in E_1$, then $z(e)=1$, so
$F(e)=-1$ and the new value is $1$. If $e$ is an edge of $Q$ not in $E_1$, then $r_f(e)\le k-3$; hence
$3\le F(e)+2\le k-1$ if $F(e)=r_f(e)$, and $3-k\le F(e)+2\le -1$ if $F(e)=r_f(e)-k$. So $F+2\chi_Q$ is a $k$-flow.
\end{proof}

\begin{proposition}[from $\Z_k$-flows to $k$-flows]\label{prop:transfer}
Let $k\ge 3$ and let $G$ be a graph. If any two $\Z_k$-flows of $G$ are joined by a chain of steps whose values are in
$\{1,-1,2,-2\}$, then $G$ is $k$-flow-connected.
\end{proposition}
\begin{proof}
Let $F$ and $F'$ be $k$-flows of $G$, and let $f$ and $f'$ be the $\Z_k$-flows formed by their residues modulo $k$. If $f=f'$, then $F$
and $F'$ are two integer lifts of $f$, and they are joined in $\cF(G,k)$ by Lemma~\ref{lem:fibre}. Otherwise let $f=f_0,f_1,\dots,f_s=f'$
with $s\ge 1$ be a chain of $\Z_k$-flows with $f_i=f_{i-1}+t_i\chi_{Q_i}$ and $t_i\in\{1,-1,2,-2\}$. If $t_i$ is
neither $1$ nor $2$, replace $Q_i$ by the same cycle with the other direction and $t_i$ by $-t_i$; so we may assume
that $t_i\in\{1,2\}$. For every $i$ from $1$ to $s$, Lemmas~\ref{lem:lift} and~\ref{lem:step1} (if $t_i=1$) or Lemma~\ref{lem:step2} (if
$t_i=2$) give an integer lift $F_{i-1}^{+}$ of $f_{i-1}$ such that $F_i^{-}=F_{i-1}^{+}+t_i\chi_{Q_i}$ is an integer
lift of $f_i$; these two $k$-flows are adjacent in $\cF(G,k)$. By Lemma~\ref{lem:fibre}, two integer lifts of the
same $\Z_k$-flow are joined in $\cF(G,k)$: $F$ is joined to $F_0^{+}$ (both are integer lifts of $f_0$), $F_i^{-}$ to
$F_i^{+}$ for $0<i<s$ (of $f_i$), and $F_s^{-}$ to $F'$ (of $f_s$). So $F$ and $F'$
are joined in $\cF(G,k)$.
\end{proof}

\begin{proof}[Proof of Theorem~\ref{thm:int}]
By Theorem~\ref{thm:cyclic} with $m=k$, any two $\Z_k$-flows of an even graph are joined by a chain of steps
with values in $\{1,-1,2,-2\}$; so Proposition~\ref{prop:transfer} shows that the graph is $k$-flow-connected.
\end{proof}

\section*{Acknowledgments}
This work was done overwhelmingly by AI models and supervised by the author. The models were Anthropic's Claude and
OpenAI's ChatGPT. The author takes full responsibility for the paper.

\end{document}